\documentclass[12pt]{article}
\usepackage{lineno,hyperref}
\usepackage{amsmath,amssymb}
\usepackage{amssymb}
\modulolinenumbers[5]

\newtheorem{Theorem}{Theorem}[section]

\newtheorem{Lemma}[Theorem]{Lemma}
\newtheorem{Corollary}[Theorem]{Corollary}
\newtheorem{Remark}[Theorem]{Remark}

\newcommand{\cvd}{\hfill$\square$ \bigskip}

\begin{document}

%
%
%
%
%
%
%
%
%

\title {Lower bounds for eigenvalues of Laplacian operator and the
clamped plate problem}

\author{\sc Zhengchao Ji \footnote{
Center of Mathematical Sciences, Zhejiang University,
             Hangzhou,  310027, People's Republic of China. E-mail: jizhengchao@zju.edu.cn} and Hongwei Xu  \footnote{
Center of Mathematical Sciences, Zhejiang University,
             Hangzhou,  310027, People's Republic of China. E-mail:xuhw@zju.edu.cn }}

\date{}
\maketitle




\begin{abstract}
In this paper, we give some lower bounds for several eigenvalues. Firstly, we investigate the eigenvalues $\lambda_i$ of the Laplace operator and prove a sharp lower bound.  Moreover, we extent this estimate of the eigenvalues to general cases.  Secondly, we study the eigenvalues $\Gamma_i$ for the clamped plate problem and deliver a sharp bound for the clamped plate problem for arbitrary dimension.
\end{abstract}
\noindent{\bf keywords:} { Laplace operator, higher eigenvalues, Weyl's asymptotic formula, P\'{o}lya conjecture,  the clamped plate problem }

\section{Introduction}\hspace{5 mm}  \setcounter{equation}{0}
Let $\Omega$ be a bounded domain with piecewise smooth boundary $\partial \Omega$ in an n-dimensional
Euclidean space $\mathbf{R}^n$. First of all, we focus on the following  Dirichlet eigenvalue problem of Laplacian
\begin{equation}\label{DLE}
\begin{cases}
\Delta u=-\lambda u  \ & \mathrm{in}\ \Omega, \\
u=0  \ &  \mathrm{on}\ \partial\Omega.
\end{cases}
\end{equation}
It is well known that the spectrum of eigenvalue problem (\ref{DLE}) is real and discrete (cf. \cite{AB,C,J,La,PPW})
\begin{alignat*}{1}
0<\lambda_1<\lambda_2\leq \lambda_3\leq \cdots\rightarrow\infty,
\end{alignat*}
where each $\lambda_i$ has finite multiplicity which is counted by its multiplicity.

Let $V(\Omega)$ be the volume of $\Omega$,  and   $\omega_n$   the volume of the unit ball in $\mathbf{R}^n$.
Then the following well-known Weyl's asymptotic formula holds
\begin{alignat*}{1}
\lambda_k\sim \frac{4\pi^2}{(\omega_nV(\Omega))^\frac{2}{n}}k^{\frac{2}{n}},\ k\rightarrow\infty,
\end{alignat*}
which implies that
\begin{alignat}{1}\label{EAF}
\frac{1}{k}\sum_{i=1}^k \lambda_i \sim \frac{n}{n+2}\frac{4\pi^2}{(\omega_nV(\Omega))^\frac{2}{n}}k^{\frac{2}{n}},\ k\rightarrow\infty.
\end{alignat}
In 1961, P\'{o}lya \cite{Po} proved that, if $n=2$ and $\Omega$ is a \textit{tiling domain} in $\mathbf{R}^2$, then
\begin{alignat*}{1}
\lambda_k \geq \frac{4\pi^2}{(\omega_nV(\Omega))^\frac{2}{n}}k^{\frac{2}{n}},\ \ \mathrm{for}\ k=1,2,\ldots,
\end{alignat*}
Based on the result above, he proposed the famous conjecture:\\

\textbf{Conjecture of P\'{o}lya} \textit{If $\Omega$ is a bounded domain in $\mathbf{R}^n$, then $k$-th eigenvalue $\lambda_k$ of the eigenvalue
problem (\ref{DLE}) satisfies}
\begin{alignat*}{1}
\lambda_k \geq \frac{4\pi^2}{(\omega_nV(\Omega))^\frac{2}{n}}k^{\frac{2}{n}},\ \ \mathrm{for}\ k=1,2,\ldots.
\end{alignat*}

During the past six decades, many mathematicians have focused on this problem and the related  topics, there are a lot of important results on this aspect (cf. \cite{Be,BE,CQ,FL,IL,KVW,KLS,Li,LP}) and we suggest  that   readers refer \cite{SY,Y} for  more details. In 1983, Li and Yau \cite{LY} verified the famous   Li-Yau inequality
\begin{alignat}{1}\label{LYE}
\frac{1}{k}\sum_{i=1}^k \lambda_i \geq \frac{n}{n+2}\frac{4\pi^2}{(\omega_nV(\Omega))^\frac{2}{n}}k^{\frac{2}{n}},\ k=1,2,\ldots.
\end{alignat}
It's seen from the asymptotic formula (\ref{EAF}),  that Li-Yau's inequality  is the best possible in the sense of the
average of eigenvalues. From (\ref{LYE}), one can derive
\begin{alignat*}{1}
\lambda_k \geq   \frac{n}{n+2}\frac{4\pi^2}{(\omega_nV(\Omega))^\frac{2}{n}}k^{\frac{2}{n}},\ \ \mathrm{for}\ k=1,2,\ldots,
\end{alignat*}
which gives a partial solution to the P\'{o}lya conjecture  with a factor $\frac{n}{n+2}$. This conjecture is still open up to now.

In \cite{M}, Melas obtained the following beautiful estimate which improves (\ref{LYE}) for $n\geq 1$ and $k\geq 1$
\begin{alignat}{1}\label{MLE}
\frac{1}{k}\sum_{i=1}^k \lambda_i \geq   \frac{n}{n+2}\frac{4\pi^2}{(\omega_nV(\Omega))^\frac{2}{n}}k^{\frac{2}{n}}+c_n\frac{V(\Omega)}{I(\Omega)},\ \ \mathrm{for}\ k=1,2,\ldots,
\end{alignat}
where $c_n$ is a positive constant depending only on  $n$ and
\begin{alignat*}{1}
I(\Omega)=\min_{a\in  \mathbf{R}^n} \int_{\Omega} |x-a|^2dx
\end{alignat*}
is called the moment of \textit{inertia} of $\Omega$. Obviously,
\begin{alignat*}{1}
I(\Omega)\geq \frac{n}{n+2}V(\Omega)\left(\frac{V(\Omega)}{\omega_n} \right)^{\frac{2}{n}}.
\end{alignat*}
Afterwards, Kova\v{r}\'{i}k, Vugalter and Weidl \cite{KVW} improved this results for $n=2$ and proved that
\begin{alignat}{1}\label{KVW}
\sum_{i=1}^k \lambda_i \geq  \frac{2\pi}{V(\Omega)}k^2+C(a_0)V(\Omega)^{-\frac{3}{2}}k^{\frac{3}{2}-\varepsilon(k)}+(1-a_0)\frac{V(\Omega)}{32I(\Omega)}k,
\end{alignat}
where $C(a_0)$ is a  positive constant depending on $a_0\in [0,1]$ and the length of the smooth parts of $\partial\Omega$,   $\varepsilon(k)=\frac{2}{\sqrt{\log_2(\frac{2\pi k}{c})}}$ and $c=\sqrt{\frac{3\pi}{14}}10^{-11}$.

The first purpose of this paper is to improve  Melas's estimate (\ref{MLE}) by giving a sharper polynomial inequality. In particular, we  improve Kova\v{r}\'{i}k, Vugalter and Weidl's inequality (\ref{KVW})  (see Theorem \ref{C1}) as follows
\begin{alignat}{1}\label{2E1}
\sum_{i=1}^k\lambda_i\geq &\frac{2\pi}{V(\Omega)}k^2+ \frac{\pi k^{\frac{3}{2}}}{ V(\Omega)}-\frac{5\pi k}{8V(\Omega)}
+\frac{\pi k^{\frac{1}{2}}}{8V(\Omega)},
\end{alignat}
for $k>k_0$, where $k_0$ is a constant. Obviously, the above inequality implies
\begin{alignat}{1}\label{2E2}
\sum_{i=1}^k\lambda_i\geq &\frac{2\pi}{V(\Omega)}k^2+ \frac{(\pi-\epsilon) k^{\frac{3}{2}}}{ V(\Omega)},
\end{alignat}
for $k\geq \left(\frac{5}{8\epsilon}\right)^2$, where $0<\epsilon<\pi$ is a constant.
 For more general cases, where $n\geq m\geq 2$ and $k\geq 1$, we obtain a lower bound for eigenvalues in Section 3, and we should mention that our result gives a sharp lower bounds by comparing Lemma \ref{KL} with the polynomial inequality in \cite{M,YY}. As a consequence of our result, we prove the Theorem \ref{Co3}, which gives that
\begin{alignat}{1}\label{mE}
\sum_{i=1}^k \lambda_i \geq &\frac{n}{n+2} \omega_n^{-\frac{2}{n}} {\alpha^{-\frac{2}{n}}}k^{\frac{2+n}{n}}+ \frac{2\omega_n^{-\frac{1}{n}}\alpha^{\frac{n-1}{n}}k^{\frac{n+1}{n}}}{(n+2)\rho}-O(k^{\frac{n-m+1}{n}}),
\end{alignat}
where $\alpha$, $\rho$ are defined by the (\ref{ar}).

The second purpose of this paper is to estimate eigenvalues of the following clamped plate problem. Let $\Omega$ be a bounded domain in $\mathbf{R}^n$. We consider the following clamped plate problem, which describes characteristic vibrations of a clamped plate:
\begin{equation*}\label{CPE}
\begin{cases}
\Delta^2 u=\Gamma u, \ & \mathrm{in}\ \Omega, \\
u=\frac{\partial u}{\partial \nu}=0, \ &  \mathrm{on}\ \partial\Omega,
\end{cases}
\end{equation*}
where $\Delta$ is the Laplacian operator and $\nu$ denotes the outward unit normal to the boundary $\partial\Omega $.
As is known,  this problem has a real and discrete spectrum (cf. \cite{A})
\begin{alignat*}{1}
0<\Gamma_1\leq\Gamma_2\leq \Gamma_3\leq \cdots\rightarrow\infty,
\end{alignat*}
where each $\Gamma_i$ has finite multiplicity which is repeated according to its multiplicity.

For the eigenvalues of the clamped plate problem, Agmon \cite{A} and Pleijel \cite{Pl} gave the
following asymptotic formula
\begin{alignat*}{1}
\Gamma_k \sim \frac{16\pi^2}{(\omega_n V(\Omega))^{\frac{4}{n}}}k^{\frac{4}{n}},\ k\rightarrow\infty.
\end{alignat*}
This implies that
\begin{alignat}{1}\label{CPAF}
\frac{1}{k}\sum_{i=1}^k\Gamma_i \sim \frac{n}{n+4} \frac{16\pi^2}{(\omega_n V(\Omega))^{\frac{4}{n}}}k^{\frac{4}{n}},\ k\rightarrow\infty.
\end{alignat}
Furthermore, Levine and Protter \cite{LP} proved that the eigenvalues of the clamped plate problem
satisfy
\begin{alignat*}{1}
\frac{1}{k}\sum_{i=1}^k\Gamma_i \geq \frac{n}{n+4} \frac{16\pi^4}{(\omega_n V(\Omega))^{\frac{4}{n}}}k^{\frac{4}{n}}.
\end{alignat*}
The formula (\ref{CPAF}) shows that the coefficient of $k^{\frac{4}{n}}$ is the best possible in the sense of the average of eigenvalues.
Later, Cheng and Wei \cite{CW1} improved the above estimate as follows:
\begin{alignat*}{1}
\frac{1}{k}\sum_{i=1}^k\Gamma_i
\geq & \frac{n}{n+4} \frac{16\pi^4}{(\omega_nV(\Omega))^{\frac{4}{n}}}k^{\frac{4}{n}}\\
&+\left( \frac{n+2}{12n(n+4)}-\frac{1}{1152n^2(n+4)}    \right)\frac{V(\Omega)}{I(\Omega)}\frac{n}{n+2}
\frac{4\pi^2}{(\omega_nV(\Omega))^{\frac{2}{n}}}k^{\frac{2}{n}}\\
&+\left( \frac{1}{576n(n+4)}-\frac{1}{27648n^2(n+2)(n+4)}    \right)\left(  \frac{V(\Omega)}{I(\Omega)}\right)^2,\\
\end{alignat*}
where $n\geq 1$ and $k\geq 1$.

Recently, by using a different method, Cheng and Wei \cite{CW2} got  better lower bounds  for eigenvalues of the clamped plate problem and proved that
\begin{alignat}{1}\label{CWE}
\begin{split}
\frac{1}{k}\sum_{i=1}^k\Gamma_i
\geq & \frac{n}{n+4} \frac{16\pi^4}{(\omega_nV(\Omega))^{\frac{4}{n}}}k^{\frac{4}{n}}+\frac{n+2}{12n(n+4)}\frac{V(\Omega)}{I(\Omega)}\frac{n}{n+2}
\frac{4\pi^2}{(\omega_nV(\Omega))^{\frac{2}{n}}}k^{\frac{2}{n}}\\
&+\frac{(n+2)^2}{1152n(n+4)^2}\left(  \frac{V(\Omega)}{I(\Omega)}\right)^2,
\end{split}
\end{alignat}
where $n\geq 2$ and $k\geq 1$.

Furthermore, they  gave  upper bounds for the sum of $\Gamma_i$,
\begin{alignat*}{1}
\frac{1}{k}\sum_{i=1}^k\Gamma_i \leq
\frac{1+\frac{4(n+4)(n^2+2n+6)}{n+2}  \frac{V(\Omega_{r_0})}{V(\Omega)}}
{\left(1-\frac{V(\Omega_{r_0})}{V(\Omega)}  \right)^{\frac{n+4}{n}}}
\frac{n}{n+4} \frac{16\pi^4}{(\omega_nV(\Omega))^{\frac{4}{n}}}k^{\frac{4}{n}},
\end{alignat*}
where $k\geq V(\Omega)r_0^n$, and
\begin{alignat*}{1}
\Omega_r=\left\{x\in\Omega \ \Big| \ \mathrm{dist}(x,\partial \Omega)<\frac{1}{r} \right\}.
\end{alignat*}

In \cite{YY1}, Yildirim and Yolcu improved Cheng and Wei's estimates by replacing the last term in the right hand of $(\ref{CWE})$ by a positive term of $k^{\frac{1}{n}}$. For any bounded open set $\Omega\subseteq R^n$, where $n\geq 2$ and $ k \geq 1$,  Yildirim and Yolcu got the following inequality.
\begin{alignat}{1}\label{TSE}
\begin{split}
\sum_{i=1}^k \Gamma_j \geq&\frac{n}{n+4}(\omega_n)^{-\frac{4}{n}}\alpha^{-\frac{4}{n}}k^{\frac{4+n}{n}}
+\frac{1}{3(n+4)}\frac{(\omega_n)^{-\frac{2}{n}}\alpha^{\frac{2n-2}{n}}k^{\frac{n+2}{n}}}{\rho^2}\\
& +\frac{2}{9(n+4)}\frac{(\omega_n)^{-\frac{1}{n}}\alpha^{\frac{3n-1}{n}}k^{\frac{n+1}{n}}}{\rho^3},
\end{split}
\end{alignat}
where
\begin{alignat}{1}\label{ar}
\alpha=\frac{V(\Omega)}{(2\pi)^n},\,\,\,\rho=\frac{V(\Omega)^{\frac{n+1}{n}}}{(2\pi)^n\omega_n^{\frac{1}{n}}}.
\end{alignat}

In Section 4, we will improve Yildirim and Yolcu's estimate (\ref{TSE}) by giving a shaper polynomial inequality. For more general cases where $n\geq m\geq 1$ and $k\geq 1$, a lower bound  will be given  as follows
\begin{alignat*}{1}
\begin{split}
\sum_{i=1}^k\Gamma_j \geq &\frac{n}{n+4} \omega_n^{-\frac{4}{n}} {\alpha^{-\frac{4}{n}}k^{1+\frac{4}{n}}}
+A_1k^{1+\frac{3}{n}}-O(k^{\frac{n-m+4}{n}}),
\end{split}
\end{alignat*}
where $A_1$ is a  positive constant depends only on $n$ and $\Omega$.

\section{Lower estimate for sums of eigenvalues}

In this section we  prove the following  theorem.

\begin{Theorem}\label{C1}
For any bounded domain $\Omega\subseteq R^n$, $n=2$ and   $ k \geq 1$ we have
\begin{alignat*}{1}
\sum_{j=1}^k\lambda_j(\Omega)\geq &\frac{n\omega_n^{-\frac{2}{n}}{\alpha}^{-\frac{2}{n}}k^{\frac{n+2}{n}}}{n+2}+ \frac{2\omega_n^{-\frac{1}{n}}\alpha^{\frac{n-1}{n}}k^{\frac{n+1}{n}}}{(n+2)\rho}\\
&-\frac{5\alpha^2k}{2(n+2)\rho^2}+\frac{\omega_n^{\frac{1}{n}}\alpha^{\frac{3n+1}{n}}k^{\frac{n-1}{n}}}{(n+2)\rho^3},
\end{alignat*}
where $\alpha$, $\rho$ are defined by (\ref{ar}).
\end{Theorem}

\begin{Remark}
Obviously, the above inequality implies that for each $\epsilon>0$ there exists a positive constant $k_\epsilon$ such that
\begin{alignat*}{1}
\sum_{i=1}^k\lambda_i\geq &\frac{2\pi}{V(\Omega)}k^2+ \frac{(\pi-\epsilon) k^{\frac{3}{2}}}{ V(\Omega)},
\end{alignat*}
for $k\geq k_\epsilon$. Hence our estimate is sharp then (\ref{KVW}).
\end{Remark}

Firstly, we  introduce some notations and definitions. For a bounded domain $\Omega$, \textit{the moment of inertia} of $\Omega$ is defined by
\begin{alignat*}{1}
I(\Omega)=\min_{a\in R^n}\int_{\Omega} |x-a|^2 dx.
\end{alignat*}
By a translation of the origin and a suitable rotation of axes, we can assume that the center
of mass is the origin and
\begin{alignat*}{1}
I(\Omega)=\int_{\Omega} |x|^2 dx.
\end{alignat*}
We now fix a $k \geq 1$ and let $u_1 ,\ldots,u_k$ denote an orthonormal set of eigenfunctions of (\ref{DLE}) corresponding to the set of eigenvalues $\lambda_1(\Omega),\ldots,\lambda_k(\Omega)$. We consider the Fourier transform of each eigenfunction
\begin{alignat*}{1}
f_j(\xi)=\hat{u}_j(\xi)=(2\pi)^{-n/2}\int_{\Omega} u_j(x)e^{ix\xi}dx.
\end{alignat*}
It seems from Plancherel's Theorem  that  $f_1 ,.\ldots, f_k$ is an orthonormal set in $R^n$. Since these eigenfunctions $u_1 ,\ldots,u_k$ are also orthonormal in $L_2(\Omega)$, Bessel's inequality implies that for every $\xi\in R^n$
\begin{alignat}{1}\label{S1}
\sum_{j=1}^k |f_j(\xi)|^2\leq (2\pi)^{-n}\int_{\Omega} |e^{ix\xi}|^2dx=(2\pi)^{-n}V(\Omega).
\end{alignat}
Since
\begin{alignat*}{1}
\nabla f_j(\xi)=(2\pi)^{-n/2}\int_{\Omega}ixu_j(x)e^{ix\xi}dx,
\end{alignat*}
we have
\begin{alignat*}{1}
\sum_{j=1}^k |\nabla f_j(\xi)|^2\leq (2\pi)^{-n/2}\int_{\Omega}|ixe^{ix\xi}|^2dx=(2\pi)^{-n}I(\Omega).
\end{alignat*}
By the boundary condition, we get
\begin{alignat*}{1}
\int_{R^n}|\xi|^2|f_j(\xi)|^2d\xi=\int_{\Omega}|\nabla u_j (x)|^2dx=\lambda_j(\Omega)
\end{alignat*}
for each $1\leq j\leq k$. Set
\begin{alignat*}{1}
F(\xi)=\sum_{j=1}^k|f_j(\xi)|^2.
\end{alignat*}
From (\ref{S1}), we have
\begin{alignat}{1}
&0\leq F(\xi)\leq (2\pi)^{-n}V(\Omega),\\
|\nabla F(\xi)|\leq 2\left( \sum_{j=1}^k|f_j(\xi)|^2 \right)^{1/2} &\left( \sum_{j=1}^k|\nabla f_j(\xi)|^2 \right)^{1/2}\leq 2(2\pi)^{-n}\sqrt{V(\Omega)I(\Omega)}\label{T1}
\end{alignat}
for each $\xi\in R^n$.  We also get
\begin{alignat}{1}
\int_{R^n}F(\xi)d\xi&=k,\\
\int_{R^n} |\xi|^2F(\xi)d\xi&=\sum_{j=1}^k \lambda_j(\Omega).\label{T2}
\end{alignat}

Assume (by approximating $F$) that the decreasing function $\phi: [0,+\infty)\rightarrow [0,(2\pi)^{-n}V(\Omega)]$  is absolutely continuous. Let $F^*(\xi) = \phi (|\xi|)$ denote the decreasing radial rearrangement of $F$.  Put $\mu(t)=|\{F^*>t \}|=|\{F>t\}|$.  It follows from the coarea formula that
\begin{alignat*}{1}
\mu(t)=\int_{t}^{(2\pi)^{-n}V(\Omega)} \int_{\{F=s \}}\frac{1}{|\nabla F|}d\sigma_sds.
\end{alignat*}
Since $F^*$ is radial, we have $\mu(\phi(s))=|\{F^*>\phi(s) \}|=\omega_n s^n$. Differentiating  both side of the above equality, we have $n\omega_n s^{n-1}={\mu}'(\phi(s))\phi'(s)$ for almost all $s$. This together with (\ref{T1}), $\rho=2(2\pi)^{-n}\sqrt{V(\Omega)I(\Omega)}$ and the isoperimetric inequality implies
\begin{alignat*}{1}
-\mu'(\phi(s)) &=\int_{\{F=\phi(s) \}} |\nabla F|^{-1} d\sigma_{\phi(s)}\\
& \geq \rho^{-1}\mathrm{Vol}_{n-1}(\{F=\phi(s)\})\\
& \geq \rho^{-1}n\omega_n s^{n-1}.
\end{alignat*}
For almost all $s$, we have
\begin{alignat}{1}\label{GI}
-\rho\leq \phi'(s)\leq 0.
\end{alignat}

Since the map $\xi \mapsto |\xi|^2$ is radial and increasing, applying (\ref{T2}), we get
\begin{alignat}{1}\label{KI}
k=\int_{R^n}F(\xi)d\xi=\int_{R^n}F^*(\xi)d\xi=n\omega_n\int_0^{\infty}s^{n-1}\phi(s)ds
\end{alignat}
and
\begin{alignat}{1}\label{LE}
\sum_{j=1}^k\lambda_j(\Omega)=\int_{R^n}|\xi|^2F(\xi)d\xi\geq \int_{R^n}|\xi|^2F^*(\xi)d\xi=n\omega_n\int_0^\infty s^{n+1}\phi(s)ds.
\end{alignat}

The following lemma will be used in the proof of  Lemma \ref{MT}.

\begin{Lemma}\label{KL}
Let $n\geq 2$, $\rho>0$, $A>0$. If $\psi: [0,+\infty)\rightarrow [0,+\infty)$ is a decreasing function (and absolutely continuous) satisfying
\begin{alignat}{1}
-\rho\leq -\psi'(s)\leq 0
\end{alignat}
and
\begin{alignat*}{1}
\int_0^\infty s^{n-1}\psi(s)ds=A.
\end{alignat*}
Then
\begin{alignat*}{1}
\int_0^{\infty} s^{n+1}\psi(s)ds \geq &\frac{(nA)^{\frac{n+2}{n}}{\psi(0)}^{-\frac{2}{n}}}{n}-\frac{s_3^3(nA)\psi(0)^2}{n(n+2)\rho^2}+\frac{s_4^4(nA)^{\frac{n-1}{n}}\psi(0)^{\frac{3n+1}{n}}}{n(n+2)\rho^3},
\end{alignat*}
where
\begin{alignat*}{1}
s_k^k=(a+1)^k-a^k\geq 1.
\end{alignat*}
\end{Lemma}

\begin{Proofp}
We choose the function $\alpha \psi(\beta t)$ for appropriate $\alpha, \beta >0$, such that $\rho = 1$ and $\psi(0) = 1$. By \cite{M} we can also assume that $B=\int_0^\infty s^{n+1}\psi(s)ds <\infty$. If we let $h(s)=-\psi^{'}(s)$ for $s\geq 0$, we have $0\leq h(s)\leq 1$ and $\int_0^\infty h(s)=\psi(0)=1.$ Moreover, integration by parts implies that
\begin{alignat*}{1}
\int_0^\infty s^{n}h(s)ds=n\int_0^\infty s^{n-1}\psi(s)ds=nA
\end{alignat*}
and
\begin{alignat*}{1}
\int_0^\infty s^{n+2}h(s)ds\leq (n+2)B.
\end{alignat*}
Next, let $0\leq a < +\infty$ satisfies that
\begin{alignat}{1}\label{DA}
\int_a^{a+1} s^{n}ds=\int_0^\infty s^{n}h(s)ds=nA.
\end{alignat}
By the same argument as in Lemma 1 of \cite{LY}, such real number $a$ exists. From \cite{M}, we have
\begin{alignat}{1}\label{BI}
(n+2)B\geq \int_0^{\infty}s^{n+2}h(s)ds\geq\int_a^{a+1}s^{n+2}ds.
\end{alignat}
To estimate the last integral we take $\tau > 0$ to be chosen later. Applying  (\ref{BI}) and integrating the both sides of the following inequality
\begin{alignat}{1}\label{K}
ns^{n+2}-(n+2)\tau^2s^n+2\tau^{n+2}\geq 2\tau^n(s-\tau)^2+4s\tau^{n-1}(s-\tau)^2, \,s\in[a,a+1],
\end{alignat}
we get
\begin{alignat*}{1}
n(n+2)B&-(n+2)\tau^2nA+2\tau^{n+2}\\
\geq& 2\tau^n \int_a^{a+1}(s-\tau)^2 +4\tau^{n-1}\int_a^{a+1}s(s-\tau)^2ds\\
\geq& 2\tau^n\left(\frac{s^3}{3}-s^2\tau+s\tau^2  \right)\bigg|_{a}^{a+1}\\
&+4\tau^{n-1}\left(\frac{s^4}{4}-\frac{2s^3\tau}{3}+\frac{s^2\tau^2}{2}  \right)\bigg|_{a}^{a+1}\\
&=2s\tau^{n+2}+2s^2\tau^{n+1}-2s^3\tau^{n}-2s^2\tau^{n+1}+s^4\tau^{n-1}\bigg|_{a}^{a+1}\\
&=2\tau^{n+2}-2s_3^3\tau^{n}+s_4^4\tau^{n-1},
\end{alignat*}
where
\begin{alignat*}{1}
s_k^k=(a+1)^k-a^k\geq 1.
\end{alignat*}
Putting, $\tau = (nA)^{1/n}$ we get
\begin{alignat*}{1}
B\geq & \frac{1}{n}(nA)^{\frac{n+2}{n}}-\frac{s_3^3}{n(n+2)}(nA)+\frac{s_4^4}{n(n+2)}(nA)^{\frac{n-1}{n}}.
\end{alignat*}
This proves Lemma \ref{KL}.

To prove (\ref{K}), we need to show that for any $\tau>0$ we have
\begin{alignat}{1}\label{KPI}
ns^{n+2}-(n+2)\tau^2s^n+2\tau^{n+2}- 2\tau^n(s-\tau)^2 -4s\tau^{n-1}(\tau-s)^2\geq 0.
\end{alignat}
Taking $t=\frac{s}{\tau}$, we define $f(t)$ (for $t>0$) by
\begin{alignat*}{1}
f(t)=nt^{n+2}-(n+2)t^n+2-2(t-1)^2-4t(t-1)^2.
\end{alignat*}
Differentiating,  $f(t)$ we have
\begin{alignat*}{1}
f'(t)& =n(n+2)t^{n+1}-(n+2)nt^{n-1}-4(t-1)-4(t-1)^2-8t(t-1)\\
&=\left[n(n+2)t^{n-2}(t+1)-12 \right]t(t-1).
\end{alignat*}
It follows from the above formula that if $n\geq 2$,  then $t=1$ is the minimum point of $f$ and   $f\geq \min\{f(1)=0,f(0)=0\}$. This implies
\begin{alignat*}{1}
f(t)\tau^{n+2}=ns^{n+2}-(n+2)\tau^2s^n+2\tau^{n+2}- 2\tau^n(s-\tau)^2 -4s\tau^{n-1}(\tau-s)^2 \geq 0.
\end{alignat*}
\cvd
\end{Proofp}

The following lemma plays important role in the proof of Theorem \ref{C1}.

\begin{Lemma}\label{MT}
For any bounded domain $\Omega\subseteq R^n$, $n\geq 2$ and  $ k \geq 1$ we have
\begin{alignat*}{1}
\sum_{j=1}^k\lambda_j(\Omega)\geq  &\omega_n^{-\frac{2}{n}}{\alpha}^{-\frac{2}{n}}k^{\frac{n+2}{n}}-\frac{s_3^3\alpha^2}{(n+2)\rho^2}k\\
&+c_1 \omega_n^{\frac{1}{n}}\frac{s_4^4\alpha^{\frac{3n+1}{n}}k^{\frac{n-1}{n}}}{(n+2)\rho^3},
\end{alignat*}
where
\begin{alignat*}{1}
c_1\leq &\min\left\{1,\max\left\{\frac{4\sqrt{2} n s_3^3k^{\frac{1}{n}}}{(3n+1)s_4^4},\frac{4\sqrt{2}(n+2)k^{\frac{3}{n}}}{(3n+1)s_4^4} \right\} \right\},\\
s_l^l=&(a+1)^{l}-a^l,
\end{alignat*}
$\alpha$, $\rho$ are defined by (\ref{ar}) and $a$ is a constant which is defined in (\ref{DA}).
\end{Lemma}

\begin{Proofp}
Applying Lemma \ref{KL} to the function $\phi$ with $A=(n\omega_n)^{-1}k, \rho=2(2\pi)^{-n}\sqrt{V(\Omega)I(\Omega)}$ and submitting it to  (\ref{LE}), we obtain
\begin{alignat}{1}\label{FI}
\begin{split}
\sum_{j=1}^k\lambda_j(\Omega)\geq  &\omega_n^{-\frac{2}{n}}{\psi(0)}^{-\frac{2}{n}}k^{\frac{n+2}{n}}-\frac{s_3^3\psi(0)^2}{(n+2)\rho^2}k\\
&+c_1 \omega_n^{\frac{1}{n}}\frac{s_4^4\psi(0)^{\frac{3n+1}{n}}k^{\frac{n-1}{n}}}{(n+2)\rho^3},
\end{split}
\end{alignat}
where $0<c_1\leq 1$ is a constant.

We observe the following facts

i)$0<\psi(0)\leq (2\pi)^{-n}V(\Omega)$,

ii)if $R$ is a positive constant such that $\omega_n R^n=V(\Omega)$, then
\begin{alignat}{1}
I(\Omega)\geq\int_{B(R)}|x|^2dx=\frac{n\omega_nR^{n+2}}{n+2}.
\end{alignat}

It follows from the above properties
\begin{alignat}{1}\label{RHO}
\rho\geq (2\pi)^{-n}\omega_n^{-\frac{1}{n}}V(\Omega)^{\frac{n+1}{n}}.
\end{alignat}
On the other hand, we consider the following function
\begin{alignat*}{1}
g(t)=g_1(t)+g_2(t),
\end{alignat*}
for $t\in(0, (2\pi)^{-n}V(\Omega)]$, where
\begin{alignat*}{1}
g_1(t)&=\omega_n^{-\frac{2}{n}}t^{-\frac{2}{n}}k^{\frac{n+2}{n}}
\end{alignat*}
and
\begin{alignat*}{1}
g_2(t)=-\frac{s_3^3t^2}{(n+2)\rho^2}k+c_1 \omega_n^{\frac{1}{n}}\frac{s_4^4t^{\frac{3n+1}{n}}k^{\frac{n-1}{n}}}{(n+2)\rho^3}.
\end{alignat*}
Then we have
\begin{alignat*}{1}
(n+2)\rho^2g_2'(t) =-{2s_3^3k t}+c_1\omega_n^{-\frac{1}{n}}\frac{s_4^4k^{\frac{n-1}{n}}}{\rho}\frac{3n+1}{n}t^{\frac{2n+1}{n}}.
\end{alignat*}

By a direct calculation, we see  from $\omega_n=\frac{2\pi^{\frac{n}{2}}}{n\Gamma(\frac{n}{2})}$ that
\begin{alignat*}{1}
\frac{\omega_n^{\frac{4}{n}}}{(2\pi)^2}\leq\frac{1}{2}.
\end{alignat*}
Therefore, in view of (\ref{RHO}), if
\begin{alignat*}{1}
c_1\leq \min\left\{1, \frac{4\sqrt{2} n s_3^3k^{\frac{1}{n}}}{(3n+1)s_4^4}  \right\},
\end{alignat*}
then $g_2(t)$ is decreasing on $(0, (2\pi)^{-n}V(\Omega)]$.
Now we consider another estimate. Seting
\begin{alignat*}{1}
G(t)=G_1(t)+G_2(t),
\end{alignat*}
where
\begin{alignat*}{1}
G_1(t)&=\omega_n^{-\frac{2}{n}}{\psi(0)}^{-\frac{2}{n}}k^{\frac{n+2}{n}}+c_1 \omega_n^{\frac{1}{n}}\frac{s_4^4\psi(0)^{\frac{3n+1}{n}}k^{\frac{n-1}{n}}}{(n+2)\rho^3}
\end{alignat*}
and
\begin{alignat*}{1}
G_2(t)=-\frac{s_3^3\psi(0)^2}{(n+2)\rho^2}k,
\end{alignat*}
we have
\begin{alignat*}{1}
G_1'(t)\rho^2 =-\frac{2}{n}\omega_n^{-\frac{2}{n}}t^{-\frac{n+2}{n}}k^{\frac{n+2}{n}}+\frac{c_1(3n+1)\omega_n^{-\frac{1}{n}}}{n}\frac{s_4^4t^{\frac{2n+1}{n}}}{(n+2)\rho^2}k^{\frac{n-1}{n}}.
\end{alignat*}
Therefore, we conclude that if
\begin{alignat*}{1}
c_1\leq \frac{4\sqrt{2}(n+2)k^{\frac{3}{n}}}{(3n+1)s_4^4},
\end{alignat*}
then $G(t)$ is decreasing on $(0, (2\pi)^{-n}V(\Omega)]$. Finally, we obtain
\begin{alignat}{1}
\begin{split}
\sum_{j=1}^k\lambda_j(\Omega)\geq  &\omega_n^{-\frac{2}{n}}{\alpha}^{-\frac{2}{n}}k^{\frac{n+2}{n}}-\frac{s_3^3\alpha^2}{(n+2)\rho^2}k\\
&+c_1 \omega_n^{\frac{1}{n}}\frac{s_4^4\alpha^{\frac{3n+1}{n}}k^{\frac{n-1}{n}}}{(n+2)\rho^3},
\end{split}
\end{alignat}
where $\alpha$, $\rho$ are defined in the (\ref{ar}) and
\begin{alignat*}{1}
c_1\leq &\min\left\{1,\max\left\{\frac{4\sqrt{2} n s_3^3k^{\frac{1}{n}}}{(3n+1)s_4^4},\frac{4\sqrt{2}(n+2)k^{\frac{3}{n}}}{(3n+1)s_4^4} \right\} \right\}.
\end{alignat*}
\cvd
\end{Proofp}

Note that $\lambda_1<\lambda_2\leq \lambda_3\leq\cdots$. This together with the above lemma implies the following estimate for higher eigenvalues.

\begin{Corollary}
For any bounded domain $\Omega\subseteq R^n$, $n\geq 2$ and any $ k \geq 1$ we have
\begin{alignat*}{1}
\lambda_k(\Omega)\geq  &\omega_n^{-\frac{2}{n}}{\alpha}^{-\frac{2}{n}}k^{\frac{2}{n}}-\frac{s_3^3\alpha^2}{(n+2)\rho^2}\\
&+c_1 \omega_n^{\frac{1}{n}}\frac{s_4^4\alpha^{\frac{3n+1}{n}}k^{\frac{-1}{n}}}{(n+2)\rho^3},
\end{alignat*}
where
\begin{alignat*}{1}
c_1\leq &\min\left\{1,\max\left\{\frac{4\sqrt{2} n s_3^3k^{\frac{1}{n}}}{(3n+1)s_4^4},\frac{4\sqrt{2}(n+2)k^{\frac{3}{n}}}{(3n+1)s_4^4} \right\} \right\},\\
s_l^l=&(a+1)^{l}-a^l,
\end{alignat*}
$\alpha$, $\rho$ are defined by (\ref{ar}).
\end{Corollary}

\begin{ProofC}\label{Re}
When $n=2$,  by applying the same argument as  in \cite{YY}, we   get that
\begin{alignat*}{1}
\int_0^{\infty} s^{n+1}\psi(s)ds \geq &\frac{(nA)^{\frac{n+2}{n}}{\psi(0)}^{-\frac{2}{n}}}{n+2}+ \frac{2(nA)^{\frac{n+1}{n}}\psi(0)^{\frac{n-1}{n}}}{n(n+2)\rho}\\
&-\frac{5(nA)\psi(0)^2}{2n(n+2)\rho^2}+\frac{(nA)^{\frac{n-1}{n}}\psi(0)^{\frac{3n+1}{n}}}{n(n+2)\rho^3}.
\end{alignat*}

Using the similar calculation in the proof of Lemma \ref{MT}, we obtain the following inequality
\begin{alignat*}{1}
\sum_{j=1}^k\lambda_j(\Omega)\geq &\frac{n\omega_n^{-\frac{2}{n}}{\alpha}^{-\frac{2}{n}}k^{\frac{n+2}{n}}}{n+2}+ \frac{2\omega_n^{-\frac{1}{n}}\alpha^{\frac{n-1}{n}}k^{\frac{n+1}{n}}}{(n+2)\rho}\\
&-\frac{5\alpha^2k}{2(n+2)\rho^2}+\frac{\omega_n^{\frac{1}{n}}\alpha^{\frac{3n+1}{n}}k^{\frac{n-1}{n}}}{(n+2)\rho^3}.
\end{alignat*}
\cvd
\end{ProofC}

\section{Lower bounds for Dirichlet eigenvalues in higher dimensions}
In this section we will give a  universal lower bound  on the sum of eigenvalues for $n\geq m+1$, where $m\geq 2$.

\begin{Theorem}\label{Co3}
When $n\geq 3$, we have the following
\begin{alignat*}{1}
\sum_{i=1}^k \lambda_i \geq &\frac{n}{n+2} \omega_n^{-\frac{2}{n}} {\alpha^{-\frac{2}{n}}}k^{\frac{2+n}{n}}+ \frac{2\omega_n^{-\frac{1}{n}}\alpha^{\frac{n-1}{n}}k^{\frac{n+1}{n}}}{(n+2)\rho}-O(k^{\frac{n-m+1}{n}}).
\end{alignat*}
\end{Theorem}

The following lemma  will be used in the proof of Lemma \ref{GMT}.
\begin{Lemma}
For an integer $n\geq m+1\geq 0$ and positive real numbers $s$ and $\tau$ we have the following
inequality:
\begin{alignat*}{1}
ns^{n+2}-(n+2)\tau^2s^n+2\tau^{n+2}-\sum_{k=1}^{m+1}2ks^{k-1}\tau^{n-k+1}(\tau-s)^2\geq 0.
\end{alignat*}
\end{Lemma}
\begin{Proofp}
Setting $t=\frac{s}{\tau}$, and putting
\begin{alignat*}{1}
f(t)=nt^{n+2}-(n+2)t^n+2-\sum_{k=1}^{m+1}2kt^{k-1}(t-1)^2,
\end{alignat*}
for $t\geq 0$,  we get
\begin{alignat}{1}
\begin{split}
f'(t)=& n(n+2)t^{n+1}-n(n+2)t^{n-1}\\
&-\left[  4(t-1)+\sum_{k=1}^m (2k(k+1)t^{k-1}(t-1)^2 + 4(k+1)t^k(t-1)     )  \right]\\
=& n(n+2)t^{n+1}-n(n+2)t^{n-1}\\
&-(t-1)\left[ 4+\sum_{k=1}^m [2k(k+1)t^{k-1}(t-1)+4(k+1)t^k  ]         \right]\\
=& n(n+2)t^{n+1}-n(n+2)t^{n-1}\\
&-(t-1)\left[ 2(m+2)(m+1)t^m +\sum_{k=1}^{m-1}( 2k(k+1)t^k-2k(k+1)t^{k-1}+ 4(k+1)t^k )   \right]\\
=& n(n+2)t^{n+1}-n(n+2)t^{n-1}-2(m+2)(m+1)t^m(t-1)  \\
=&t^m(t-1)\left[ n(n+2)t^{n-m-1}(t+1)-2(m+2)(m+1)  \right].
\end{split}
\end{alignat}
It follows from the above formula that if $n \geq m+1$, then $t = 1$ is the minimum point of $f(t)$ and
 $f \geq \min\{f(1) = 0,f(0) = 0\}$. So,  we get
\begin{alignat*}{1}
\tau^{n+2}f(t)=ns^{n+2}-(n+2)\tau^2s^n-\sum_{k=1}^{m+1}2ks^{k-1}\tau^{n-k+1}(\tau-s)^2\geq 0.
\end{alignat*}
\cvd
\end{Proofp}

The following lemma plays important role in the proof of  Theorem \ref{Co3}.

\begin{Lemma}\label{GMT}
For any bounded domain $\Omega\subseteq R^n$, $n\geq m+1\geq 3$ and  $ k \geq 1$, we have
\begin{alignat*}{1}
\sum_{i=1}^k \lambda_i \geq & \omega_n^{-\frac{2}{n}} {\beta^{-\frac{2}{n}}}k^{\frac{n+2}{n}}-\frac{2\omega_n^{\frac{m-1}{n}}S_{m+2}\beta^{\frac{(m+1)n+m-1}{n}}}{(n+2)\rho^{m+1}}k^{\frac{n-m+1}{n}}\\
&+c_2\frac{2\omega_n^{\frac{m}{n}}(m+1)S_{m+3}\beta^{\frac{(m+2)n+m}{n}}}{(n+2)(m+3)\rho^{m+2}}k^{\frac{n-m}{n}},
\end{alignat*}
where
\begin{alignat*}{1}
\begin{split}
c_2\leq & \min \left\{1,  \frac{(m+1)n+m-1 }{(m+2)n+m}\frac{\sqrt{2}S_{m+2}}{S_{m+3}}\frac{m+3}{m+1}k^{\frac{1}{n}} \right\},\\
S_{k}=&(a+1)^k-a^k,
\end{split}
\end{alignat*}
$\beta=\frac{V(\Omega)}{(2\pi)^n},\,\,\,\,\rho=\frac{V(\Omega)^{\frac{n+1}{n}}}{(2\pi)^n\omega_n^{\frac{1}{n}}}$ and $a$ is a constant defined in (\ref{DA}).
\end{Lemma}

\begin{Proofp}For $l\geq 0$, $\tau\geq \frac{1}{2}$ and $a \geq 0$, we have
\begin{alignat}{1}
\begin{split}
\int_a^{a+1}s^l(\tau-s)^2ds &=\frac{s^{l+3}}{l+3}-\frac{2s^{l+2}}{l+2}\tau+\frac{s^{l+1}}{l+1}\tau^{2}\bigg|_a^{a+1}\\
&=\frac{S_{l+3}}{l+3}-\frac{2S_{l+2}}{l+2}\tau+\frac{S_{l+1}}{l+1}\tau^{2},
\end{split}
\end{alignat}
where
\begin{alignat*}{1}
S_j=(a+1)^j-a^j\geq 1.
\end{alignat*}

Therefore, we get
\begin{alignat*}{1}
n(n+2)B-(n+2)\tau^2nA+2\tau^{n+2} \geq &  \sum_{k=1}^{m+1}2k\tau^{n-k+1}\left(\frac{S_{k+2}}{k+2}-\frac{2S_{k+1}}{k+1}\tau+\frac{S_{k}}{k}\tau^{2}\right).
\end{alignat*}
From
\begin{alignat*}{1}
\sum_{k=1}^{m+1}&2k\tau^{n-k+1}\left(\frac{S_{k}}{k}\tau^{2}-\frac{2S_{k+1}}{k+1}\tau +\frac{S_{k+2}}{k+2}\right)\\
=&2\tau^{n+2}+2\sum_{k=1}^{m}S_{k+1}\tau^{n-k+2}-2\sum_{k=1}^{m+1}\frac{2kS_{k+1}}{k+1}\tau^{n-k+2}+2\sum_{k=1}^{m+1}\frac{kS_{k+2}}{k+2}\tau^{n-k+1}\\
=&2\tau^{n+2}+2S_2\tau^{n+1}+\frac{2mS_{m+2}}{m+2}\tau^{n-m+1}+\frac{2(m+1)S_{m+3}}{m+3}\tau^{n-m}\\&
-2S_2\tau^{n+1}- \frac{4(m+1)S_{m+2}}{m+2}\tau^{n-m+1}\\
&+2\sum_{k=2}^{m }\left(1+\frac{k-1}{k+1}-\frac{2k}{k+1} \right)S_{k+1}\tau^{n-k+2}\\
=&2\tau^{n+2}+2S_2\tau^{n+1}+\frac{2mS_{m+2}}{m+2}\tau^{n-m+1}+\frac{2(m+1)S_{m+3}}{m+3}\tau^{n-m}\\
&-2S_2\tau^{n+1}- \frac{4(m+1)S_{m+2}}{m+2}\tau^{n-m+1}\\
=&2\tau^{n+2}-2S_{m+2}\tau^{n-m+1}+\frac{2(m+1)S_{m+3}}{m+3}\tau^{n-m},
\end{alignat*}
and
\begin{alignat*}{1}
\sum_{k=2}^{m }\left(1+\frac{k-1}{k+1}-\frac{2k}{k+1} \right)S_{k+1}\tau^{n-k+2}=0,
\end{alignat*}
we obtain
\begin{alignat*}{1}
n(n+2)B-(n+2)\tau^2nA+2\tau^{n+2} \geq &2\tau^{n+2}-2S_{m+2}\tau^{n-m+1}\\
&+\frac{2(m+1)S_{m+3}}{m+3}\tau^{n-m}.
\end{alignat*}

Choosing $\tau=(nA)^{\frac{1}{n}}$, we get
\begin{alignat}{1}\label{FKI}
\begin{split}
B\geq& \frac{(nA)^{\frac{n+2}{n}}}{n}-\frac{2S_{m+2}(nA)^{\frac{n-m+1}{n}}}{n(n+2)}+\frac{2(m+1)S_{m+3}(nA)^{\frac{n-m}{n}}}{n(n+2)(m+3)}.
\end{split}
\end{alignat}
It follows from (\ref{FKI}) that
\begin{alignat}{1}
\begin{split}
\int_0^\infty  s^{n+1}\psi(s)ds\geq & \frac{(nA)^{\frac{n+2}{n}}\psi(0)^{-\frac{2}{n}}}{n}-\frac{2S_{m+2}(nA)^{\frac{n-m+1}{n}}\psi(0)^{\frac{(m+1)n+m-1}{n}}}{n(n+2)\rho^{m+1}}\\
&+\frac{2(m+1)S_{m+3}(nA)^{\frac{n-m}{n}}\psi(0)^{\frac{(m+2)n+m}{n}}}{n(n+2)(m+3)\rho^{m+2}}.
\end{split}
\end{alignat}
From (\ref{LE}), we know
\begin{alignat*}{1}
\begin{split}
\sum_{i=1}^k \lambda_i \geq & n\omega_n \int_0^\infty  s^{n+1}\psi(s)ds\\
\geq &\omega_n {(nA)^{\frac{n+2}{n}}\psi(0)^{-\frac{2}{n}}}-\frac{2\omega_nS_{m+2}(nA)^{\frac{n-m+1}{n}}\psi(0)^{\frac{(m+1)n+m-1}{n}}}{(n+2)\rho^{m+1}}\\
&+\frac{2\omega_n(m+1)S_{m+3}(nA)^{\frac{n-m}{n}}\psi(0)^{\frac{(m+2)n+m}{n}}}{(n+2)(m+3)\rho^{m+2}}.
\end{split}
\end{alignat*}
In view of $A=\frac{k}{n\omega_n}$, we have
\begin{alignat}{1}
\begin{split}
\sum_{i=1}^k \lambda_i \geq & \omega_n^{-\frac{2}{n}} {\psi(0)^{-\frac{2}{n}}}k^{\frac{n+2}{n}}-\frac{2\omega_n^{\frac{m-1}{n}}S_{m+2}\psi(0)^{\frac{(m+1)n+m-1}{n}}}{(n+2)\rho^{m+1}}k^{\frac{n-m+1}{n}}\\
&+c_2\frac{2\omega_n^{\frac{m}{n}}(m+1)S_{m+3}\psi(0)^{\frac{(m+2)n+m}{n}}}{(n+2)(m+3)\rho^{m+2}}k^{\frac{n-m}{n}},
\end{split}
\end{alignat}
where $0<c_2\leq 1$ is a constant.

When $m=1$, we complete the proof of Theorem \ref{GMT} in section 2. We assume that $m\geq 2$. Putting
\begin{alignat*}{1}
g(t)=g_1(t)+g_2(t),
\end{alignat*}
where
\begin{alignat*}{1}
g_1(t)&= \omega_n^{-\frac{2}{n}} {t^{-\frac{2}{n}}}k^{\frac{n+2}{n}}
\end{alignat*}
and
\begin{alignat*}{1}
g_2(t)=&-\frac{2\omega_n^{\frac{m-1}{n}}S_{m+2}t^{\frac{(m+1)n+m-1}{n}}}{(n+2)\rho^{m+1}}k^{\frac{n-m+1}{n}}\\
&+c_2\frac{2\omega_n^{\frac{m}{n}}(m+1)S_{m+3}t^{\frac{(m+2)n+m}{n}}}{(n+2)(m+3)\rho^{m+2}}k^{\frac{n-m}{n}},
\end{alignat*}
we have
\begin{alignat*}{1}
\frac{(n+2)\rho^{m+1}\omega_n^{\frac{m}{n}}g_2'(t)}{2k^{\frac{n-m}{n}}}=&-\frac{(m+1)n+m-1}{n}\omega_n^{-\frac{1}{n}}S_{m+2}t^{\frac{mn+m-1}{n}}k^{\frac{1}{n}}\\
&+c_2\frac{(m+2)n+m}{n}    \frac{(m+1)S_{m+3}}{(m+3)\rho}t^{\frac{(m+1)n+m}{n}  }.
\end{alignat*}
When
\begin{alignat}{1}
\begin{split}
c_2\leq \frac{(m+1)n+m-1 }{(m+2)n+m}\frac{\sqrt{2}S_{m+2}}{S_{m+3}}\frac{m+3}{m+1}k^{\frac{1}{n}},
\end{split}
\end{alignat}
we get that $g_2(t)$ is decreasing on $(0, (2\pi)^{-n}V(\Omega)]$ by using the following formulas
\begin{alignat*}{1}
nA&=\frac{k}{\omega_n},\\
\rho\geq (2\pi)^{-n}&\omega_n^{-\frac{1}{n}}V(\Omega)^{\frac{n+1}{n}}.\\
\end{alignat*}
Hence  $g(t)$ is also decreasing on $(0, (2\pi)^{-n}V(\Omega)]$.
This implies
\begin{alignat}{1}
\begin{split}
\sum_{i=1}^k \lambda_i \geq & \omega_n^{-\frac{2}{n}} {\psi(0)^{-\frac{2}{n}}}k^{\frac{n+2}{n}}-\frac{2\omega_n^{\frac{m-1}{n}}S_{m+2}\psi(0)^{\frac{(m+1)n+m-1}{n}}}{(n+2)\rho^{m+1}}k^{\frac{n-m+1}{n}}\\
&+c_2\frac{2\omega_n^{\frac{m}{n}}(m+1)S_{m+3}\psi(0)^{\frac{(m+2)n+m}{n}}}{(n+2)(m+3)\rho^{m+2}}k^{\frac{n-m}{n}},
\end{split}
\end{alignat}
where
\begin{alignat*}{1}
\psi(0)&=\frac{V(\Omega)}{(2\pi)^n},\\
\end{alignat*}
and
\begin{alignat*}{1}
\rho&=\frac{V(\Omega)^{\frac{n+1}{n}}}{(2\pi)^n\omega_n^{\frac{1}{n}}}.
\end{alignat*}

\cvd
\end{Proofp}

From the above lemma, we have the following universal lower bounds for higher eigenvalues.

\begin{Corollary}
For any bounded domain $\Omega\subseteq R^n$, $n\geq m+1\geq 3$ and $ k \geq 1$ we have
\begin{alignat*}{1}
\lambda_k \geq & \omega_n^{-\frac{2}{n}} {\beta^{-\frac{2}{n}}}k^{\frac{2}{n}}-\frac{2\omega_n^{\frac{m-1}{n}}S_{m+2}\beta^{\frac{(m+1)n+m-1}{n}}}{(n+2)\rho^{m+1}}k^{\frac{-m+1}{n}}\\
&+c_2\frac{2\omega_n^{\frac{m}{n}}(m+1)S_{m+3}\beta^{\frac{(m+2)n+m}{n}}}{(n+2)(m+3)\rho^{m+2}}k^{\frac{-m}{n}},
\end{alignat*}
where
\begin{alignat*}{1}
\begin{split}
c_2\leq & \min \left\{1,  \frac{(m+1)n+m-1 }{(m+2)n+m}\frac{\sqrt{2}S_{m+2}}{S_{m+3}}\frac{m+3}{m+1}k^{\frac{1}{n}} \right\},\\
S_{k}=&(a+1)^k-a^k,
\end{split}
\end{alignat*}
$\beta=\frac{V(\Omega)}{(2\pi)^n},\,\,\,\,\rho=\frac{V(\Omega)^{\frac{n+1}{n}}}{(2\pi)^n\omega_n^{\frac{1}{n}}}$ and $a$ is a constant defined in (\ref{DA}).
\end{Corollary}

\begin{ProofCo}
Applying Lemma \ref{GMT} and using the similar strategy as in the proof of  Theorem \ref{C1}, we prove the desired inequality.

\cvd
\end{ProofCo}

\section{A universal lower bound on  eigenvalues of the clamped plate problem}
In this section, let $a$ be a real number defined  as in (\ref{DA}) where we replace $n$ by $n+3$. We will give a  universal lower bounds on the sum of eigenvalues for $n\geq m$, where $m\geq 1$.

\begin{Theorem}\label{CPT}
For any bounded domain $\Omega\subseteq R^n$, $n\geq m\geq 1$ and   $ k \geq 1$ we have\\
\begin{alignat}{1}
\begin{split}
\sum_{i=1}^k\Gamma_i \geq &\frac{n}{n+4} \omega_n^{-\frac{4}{n}} {\alpha^{-\frac{4}{n}}k^{1+\frac{4}{n}}}
+A_1k^{1+\frac{3}{n}}-O(k^{\frac{n-m+4}{n}}),
\end{split}
\end{alignat}
where $A_1>0$ is a constant only depends on $n$ and $\Omega$.
\end{Theorem}

Next, we recall the definition and serval properties of the symmetric decreasing rearrangements.  Let $\Omega\subset R^n $ be a bounded domain. Its symmetric rearrangement $\Omega^*$ is the open ball with the same volume as $\Omega$,
\begin{alignat*}{1}
\Omega^*=\left\{x\in R^n| \,\, |x|<\left(\frac{V(\Omega)}{\omega_n}  \right)   \right\}.
\end{alignat*}
By using a symmetric rearrangement of $\Omega$, we have
\begin{alignat}{1}\label{II}
I(\Omega)=\int_{\Omega} |x|^2 dx\geq \int_{\Omega^*} |x|^2 dx=\frac{n}{n+2}V(\Omega)\left(\frac{V(\Omega)}{\omega_n} \right)^{\frac{2}{n}}.
\end{alignat}
Then we have
\begin{alignat}{1}\label{KE}
\int_{R^n}|x|^4F(x)dx\geq \int_{R^n}|x|^4F^{*}(x)dx=n\omega_n\int_0^{\infty} s^{n+3}\phi(s)ds.
\end{alignat}

The following lemma s useful in the proof of Lemma \ref{TLT}.
\begin{Lemma}\label{SKL}
For   integers $n\geq m\geq 1$ and positive real numbers $s$ and $\tau$, we have the following
inequality:
\begin{alignat}{1}\label{LGIE}
ns^{n+4}-(n+4)\tau^4s^n+4\tau^{n+4}-\sum_{k=1}^{m}4ks^{k-1}\tau^{n-k+3}(\tau-s)^2\geq 0.
\end{alignat}
\end{Lemma}
\begin{Proofp}
Taking $t=\frac{s}{\tau}$, and putting $f(t)$
\begin{alignat*}{1}
f(t)=nt^{n+4}-(n+4)t^n+4-4(t-1)^2 -  \sum_{k=2}^{m}4kt^{k-1}(t-1)^2,
\end{alignat*}
for $t\geq 0$,  we get
\begin{alignat}{1}
\begin{split}
f'(t)=&  n(n+4)t^{n+3}-n(n+4)t^{n-1}\\
&  -\left[8(t-1)+ \sum_{k=2}^{m}4k(k-1)t^{k-2}(t-1)^2+\sum_{k=2}^{m}8kt^{k-1}(t-1)   \right] \\
=&  n(n+4)t^{n+3}-n(n+4)t^{n-1}\\
&  -(t-1)\left[8+ \sum_{k=2}^{m}4k(k-1)t^{k-2}(t-1)+\sum_{k=2}^{m}8kt^{k-1}   \right] \\
=&  n(n+4)t^{n+3}-n(n+4)t^{n-1}\\
&  -(t-1)\left[8+ \sum_{k=2}^{m}4k(k-1)t^{k-1} -\sum_{k=2}^{m}4k(k-1)t^{k-2}    +\sum_{k=2}^{m}8kt^{k-1}   \right] \\
=&  n(n+4)t^{n+3}-n(n+4)t^{n-1}\\
&  -(t-1)\left[4m(m+1)t^{m-1} + \sum_{k=3}^{m}(4(k-1)(k-2) -4k(k-1)+8(k-1) )t^{k-2}  \right] \\
=&n(n+4)t^{n+3}-n(n+4)t^{n-1}-4m(m+1)t^{m-1}(t-1)\\
=&\left[ n(n+4)t^{n-m}(t^2+1)(t+1)-4m(m+1)  \right]t^{n-m}(t-1).
\end{split}
\end{alignat}
From the above formula, it is clear that when $n \geq m$, we have $t = 1$ is the minimum point of $f(t)$ and
then $f \geq \min\{f(1) = 0,f(0) = 0\}$. We get
\begin{alignat*}{1}
\tau^{n+4}f(t)=ns^{n+4}-(n+4)\tau^4s^n-\sum_{k=1}^{m}4ks^{k-1}\tau^{n-k+3}(\tau-s)^2\geq 0.
\end{alignat*}
\cvd
\end{Proofp}

The following lemma palys important role in the proof of Theorem \ref{CPT}.

\begin{Lemma}\label{TLT}
For any bounded domain $\Omega\subseteq R^n$, $n\geq m\geq 1$ and   $ k \geq 1$ we have\\

(1) When $n=1$ and
\begin{alignat*}{1}
\frac{2\sqrt{2}S_3}{5} \leq k,
\end{alignat*}
we have
\begin{alignat}{1}
\begin{split}
\sum_{i=1}^k\Gamma_i \geq & \omega_n^{-\frac{4}{n}} {\alpha^{-\frac{4}{n}}k^{1+\frac{4}{n}}}
-\omega_n^{\frac{m-4}{n}} \frac{4S_{m+2}}{(n+4)\rho^{m}}\alpha^{\frac{mn+m-4}{n}}k^{\frac{n-m+4}{n}}\\
&+\omega_n^{\frac{m-3}{n}} \frac{4mS_{m+2}}{(n+4)(m+2)\rho^{m+1}}\alpha^{\frac{(m+1)n+m-3}{n}}k^{\frac{n-m+3}{n}},
\end{split}
\end{alignat}\\
where $\alpha$, $\rho$ are defined by (\ref{ar}) and
\begin{alignat*}{1}
S_l=(a+1)^l-a^l.
\end{alignat*}
\\
(2) When $m\geq 2$, we have
\begin{alignat}{1}
\begin{split}
\sum_{i=1}^k\Gamma_i \geq & \omega_n^{-\frac{4}{n}} {\alpha^{-\frac{4}{n}}k^{1+\frac{4}{n}}}
-\omega_n^{\frac{m-4}{n}} \frac{4S_{m+2}}{(n+4)\rho^{m}}\alpha^{\frac{mn+m-4}{n}}k^{\frac{n-m+4}{n}}\\
&+c_3\omega_n^{\frac{m-3}{n}} \frac{4mS_{m+2}}{(n+4)(m+2)\rho^{m+1}}\alpha^{\frac{(m+1)n+m-3}{n}}k^{\frac{n-m+3}{n}},
\end{split}
\end{alignat}
where
\begin{alignat*}{1}
c_3\leq \min \left\{ 1, \frac{2^{\frac{m+1}{2}}(n+2)(m+2)}{S_{m+2}[(m+1)n+m-3]}k^{\frac{m+1}{n}} \right\}.
\end{alignat*}
\end{Lemma}

\begin{Proofp}Let $\{u_j\}_{j=1}^{\infty}$ be  the eigenfunction corresponding to the eigenvalue $\Gamma_j$, $j=1,2.\ldots$  which satisfy
\begin{equation*}
        \begin{cases}
           \Delta^2 u_j=\Gamma_j u_j,\,\,&\mathrm{in}\,\,\Omega,\\
           u_j=\frac{\partial u_j}{\partial \nu}=0,\,\,&\mathrm{on}\,\,\partial\Omega,\\
           \int_{\Omega}u_i(x) u_j(x)dx=\delta_{ij},\,\,&\mathrm{for}\,\,\mathrm{any}\,\,i,j.
        \end{cases}
\end{equation*}
Thus, $\{ u_j \}_{j=1}^{\infty}$ forms an orthonormal basis of $L^2(\Omega)$. We define a function $\varphi_j$ by
\begin{equation*}
\varphi_j(x)=\begin{cases}
            u_j(x), \ &x \in \Omega, \\
            0, \ &x \in\mathbf{R}^n \backslash\Omega.
        \end{cases}
\end{equation*}
Denote by $\widehat{\varphi}_j(z)$ the Fourier transform of $\varphi_j (x)$. For any $z \in \mathbf{R}^n$, we have
\begin{alignat*}{1}
\widehat{\varphi}_j(z)=(2\pi)^{-\frac{n}{2}}\int_{\mathbf{R}^n}\varphi_j(x)e^{i\langle x,z\rangle}dx=(2\pi)^{-\frac{n}{2}}\int_{\Omega}u_j(x)e^{i\langle x,z\rangle}dx.
\end{alignat*}
By the Plancherel formula, we have
\begin{alignat}{1}\label{VI}
\int_{\mathbf{R}^n}\widehat{\varphi}_i(z)\widehat{\varphi}_j(z)=\delta_{ij}
\end{alignat}
for any $i,j$. Since $\{ u_j \}_{j=1}^{\infty}$ is an orthonormal basis in $L^2(\Omega)$, the Bessel inequality implies that
\begin{alignat*}{1}
\sum_{j=1}^k |\widehat{\varphi}_j(z) |^2\leq(2\pi)^{-n}\int_{\Omega}|e^{i\langle x,z \rangle} |^2dx=(2\pi)^{-n}V(\Omega).
\end{alignat*}
For each  $j = 1,\ldots,k$, we deduce from the divergence theorem and $u_j|_{\partial \Omega}=\frac{\partial u_j}{\partial \nu}|_{\partial \Omega}=0$ that
\begin{alignat*}{1}
z_p^2\widehat{\varphi}_j(z)&=(2\pi)^{-\frac{n}{2}}\int_{\mathbf{R}^n}\varphi_j(x)(-i)^2 \frac{\partial^2e^{i\langle x,z \rangle} }{\partial x_p^2}dx\\
&=-(2\pi)^{-\frac{n}{2}}\int_{\mathbf{R}^n}\frac{\partial^2\varphi_j(x) }{\partial x_p^2} e^{i\langle x,z \rangle}dx\\
&=-\frac{\widehat{\partial^2\varphi_j}}{\partial x_p^2}(z).
\end{alignat*}
It follows from the Parseval's identity that
\begin{alignat}{1}\label{GE}
\begin{split}
\int_{\mathbf{R}^n}|z|^4|\widehat{\varphi}_j(z) |^2dz &=\int_{\mathbf{R}^n}(|z|^2|\widehat{\varphi}_j(z) |)^2dz\\
&=\int_{\Omega}|\Delta u_j(x)|^2dx\\
&=\Gamma_j.
\end{split}
\end{alignat}
Since
\begin{alignat*}{1}
\nabla \widehat{\varphi}_j(z) =(2\pi)^{-\frac{n}{2}}\int_{\Omega} ixu_j(x)e^{i\langle x,z \rangle}dx,
\end{alignat*}
we obtain
\begin{alignat}{1}\label{NV}
\sum_{j=1}^k  |\nabla \widehat{\varphi}_j(z) |^2\leq (2\pi)^{-n}\int_{\Omega} |ixe^{i\langle x,z \rangle}|^2dx=(2\pi)^{-n}I(\Omega).
\end{alignat}
Putting
\begin{alignat*}{1}
h(z):=\sum_{j=1}^k  | \widehat{\varphi}_j(z) |^2,
\end{alignat*}
one derives from (\ref{VI}) that $0 \leq h(z) \leq (2\pi)^{-n}V(\Omega)$. It follows from (\ref{NV}) and the Cauchy-Schwarz inequality that
\begin{alignat*}{1}
|\nabla h(z)| & \leq 2\left(\sum_{j=1}^k  | \widehat{\varphi}_j(z) |^2   \right)^{\frac{1}{2}}\left(\sum_{j=1}^k  |\nabla \widehat{\varphi}_j(z) |^2   \right)^{\frac{1}{2}}\\
& \leq 2(2\pi)^{-n}\sqrt{V(\Omega)I(\Omega)}
\end{alignat*}
for every $z\in \mathbf{R }^n$. From the Parseval's identity, we derive
\begin{alignat}{1}\label{KII}
\int_{\mathbf{R }^n}h(z)dz=\sum_{j=1}^k  \int_{\Omega}|u_j(x)|^2dx=k.
\end{alignat}
Applying the symmetric decreasing rearrangement to $h(z)$ and noting that $\zeta=\sup |\nabla h|\leq  2(2\pi)^{-n}\sqrt{V(\Omega)I(\Omega)}:=\eta,$  we see from (\ref{GI})
\begin{alignat*}{1}
-\eta\leq -\zeta\leq\phi{'}(s)\leq 0
\end{alignat*}
for almost every $s$. According to (\ref{KE}) and (\ref{GE}), we infer
\begin{alignat}{1}\label{KEE}
\begin{split}
\sum_{i=1}^k\Gamma_j&=\int_{\mathbf{R}^n}|z|^4h(z)dz\\
&\geq \int_{\mathbf{R}^n}|z|^4h^{*}(z)dz\\
&=n\omega_n \int_{0}^{\infty} s^{n+3}\phi(s)ds.
\end{split}
\end{alignat}
In order to apply Lemma \ref{SKL}, from (\ref{KE}) and the definition of $A$, we take
\begin{alignat}{1}\label{TE}
\psi(s)=\phi(s),\,\, A=\frac{k}{n\omega_n},\,\,\eta=2(2\pi)^{-n}\sqrt{V(\Omega)I(\Omega)}.
\end{alignat}
From (\ref{II}), we deduce that
\begin{alignat}{1}\label{EI}
\rho\geq 2(2\pi)^{-n}\left(\frac{n}{n+2}  \right)^{\frac{1}{2}}\omega_n^{-\frac{1}{n}}V(\Omega)^{\frac{n+1}{n}}.
\end{alignat}
On the other hand, $0 < \phi(0) \leq \sup h^*(z) = \sup h(z) \leq (2\pi)^{-n} V(\Omega)$.

For any $k\geq 1$ and $a \geq 0$, we have
\begin{alignat}{1}
\begin{split}
\int_a^{a+1}s^{k-1}(\tau-s)^2ds &=\frac{s^{k+2}}{k+2}-\frac{2s^{k+1}}{k+1}\tau+\frac{s^{k}}{k}\tau^{2}\bigg|_a^{a+1}\\
&=\frac{S_{k+2}}{k+2}-\frac{2S_{k+1}}{k+1}\tau+\frac{S_{k}}{k}\tau^{2},
\end{split}
\end{alignat}
where
\begin{alignat*}{1}
S_l=(a+1)^l-a^l.
\end{alignat*}
From the above lemma, integrating the both sides of (\ref{LGIE}) over $[a,a+1]$,  we get
\begin{alignat*}{1}
n(n+4)D'-(n+4)\tau^4 nA+4\tau^{n+4}\geq \sum_{k=1}^m4k\tau^{n-k+3}\left(\frac{S_{k}}{k}\tau^{2}-\frac{2S_{k+1}}{k+1}\tau+\frac{S_{k+2}}{k+2}\right).
\end{alignat*}
From
\begin{alignat*}{1}
\sum_{k=1}^{m}&4k\tau^{n-k+3}\left(\frac{S_{k}}{k}\tau^{2}-\frac{2S_{k+1}}{k+1}\tau +\frac{S_{k+2}}{k+2}\right)\\
=&4\tau^{n+4}+4\sum_{k=1}^{m-1}S_{k+1}\tau^{n-k+4}-4\sum_{k=1}^{m}\frac{2kS_{k+1}}{k+1}\tau^{n-k+4}+4\sum_{k=1}^{m }\frac{kS_{k+2}}{k+2}\tau^{n-k+3}\\
=&4\tau^{n+4}+4S_2\tau^{n+3}+\frac{4mS_{m+2}}{m+2}\tau^{n-m+3}+\frac{4(m-1)S_{m+1}}{m+1}\tau^{n-m+4}\\&
-4S_2\tau^{n+3}- \frac{8mS_{m+2}}{m+1}\tau^{n-m+4}\\
&+4\sum_{k=2}^{m-1 }\left(1+\frac{k-1}{k+1}-\frac{2k}{k+1} \right)S_{k+1}\tau^{n-k+4}\\
=&4\tau^{n+4}+4S_2\tau^{n+3}+\frac{4mS_{m+2}}{m+2}\tau^{n-m+3}+\frac{4(m-1)S_{m+1}}{m+1}\tau^{n-m+4}\\
&-4S_2\tau^{n+3}- \frac{8mS_{m+2}}{m+1}\tau^{n-m+4}\\
=&4\tau^{n+4}-4S_{m+2}\tau^{n-m+4}+\frac{4mS_{m+2}}{m+2}\tau^{n-m+3},
\end{alignat*}
and
\begin{alignat*}{1}
4\sum_{k=2}^{m-1 }\left(1+\frac{k-1}{k+1}-\frac{2k}{k+1} \right)S_{k+1}\tau^{n-k+4}=0,
\end{alignat*}
we get
\begin{alignat*}{1}
n(n+4)D'-(n+4)\tau^4 nA+4\tau^{n+4}\geq& 4\tau^{n+4}-4S_{m+2}\tau^{n-m+4}\\
&+\frac{4mS_{m+2}}{m+2}\tau^{n-m+3}.
\end{alignat*}
This implies that
\begin{alignat*}{1}
n(n+4)D'\geq &(n+4)\tau^4(nA) -4S_{m+2}\tau^{n-m+4}\\
&+\frac{4mS_{m+2}}{m+2}\tau^{n-m+3}.
\end{alignat*}

Taking $\tau=(nA)^{\frac{1}{n}}$, we get
\begin{alignat*}{1}
D'\geq& \frac{(nA)}{n}\tau^4 -\frac{4S_{m+2}}{n(n+4)}\tau^{n-m+4}\\
&+\frac{4mS_{m+2}}{n(n+4)(m+2)}\tau^{n-m+3}\\
\geq&\frac{(nA)^{\frac{n+4}{n}}}{n}-\frac{4S_{m+2}(nA)^{\frac{n-m+4}{n}}}{n(n+4)}\\
&+\frac{4mS_{m+2}(nA)^{\frac{n-m+3}{n}}}{n(n+4)(m+2)}.
\end{alignat*}
Then, we get
\begin{alignat}{1}
\begin{split}
\int_0^{\infty}s^{n+3}\psi(s)ds\geq& \frac{(nA)^{1+\frac{4}{n}}}{n}\psi(0)^{-\frac{4}{n}}-\frac{4S_{m+2}(nA)^{\frac{n-m+4}{n}}}{n(n+4)\rho^{m}}\psi(0)^{\frac{mn+m-4}{n}}\\
&+\frac{4mS_{m+2}(nA)^{\frac{n-m+3}{n}}}{n(n+4)(m+2)\rho^{m+1}}\rho^{\frac{(m+1)n+m-3}{n}}.
\end{split}
\end{alignat}
According to (\ref{KE}), (\ref{GE}) and the above inequality, we conclude
\begin{alignat}{1}
\begin{split}
\sum_{i=1}^k\Gamma_j =&\int_{\mathbf{R}^n}|z|^4h(z)dz\\
\geq& \int_{\mathbf{R}^n}|z|^4h^{*}(z)dz\\
=&n\omega_n \int_{0}^{\infty} s^{n+3}\phi(s)ds\\
\geq& n\omega_n  \frac{(nA)^{1+\frac{4}{n}}}{n}\psi(0)^{-\frac{4}{n}}-n\omega_n \frac{4S_{m+2}(nA)^{\frac{n-m+4}{n}}}{n(n+4)\rho^{m}}\psi(0)^{\frac{mn+m-4}{n}}\\
&+n\omega_n \frac{4mS_{m+2}(nA)^{\frac{n-m+3}{n}}}{n(n+4)(m+2)\rho^{m+1}}\rho^{\frac{(m+1)n+m-3}{n}}\\
=& \omega_n {(nA)^{1+\frac{4}{n}}}\psi(0)^{-\frac{4}{n}}-\omega_n \frac{4S_{m+2}(nA)^{\frac{n-m+4}{n}}}{(n+4)\rho^{m}}\psi(0)^{\frac{mn+m-4}{n}}\\
&+\omega_n \frac{4mS_{m+2}(nA)^{\frac{n-m+3}{n}}}{(n+4)(m+2)\rho^{m+1}}\psi(0)^{\frac{(m+1)n+m-3}{n}}.
\end{split}
\end{alignat}
For $m=1$ and $n=1$, we define $f(t)$ as follows
\begin{alignat*}{1}
f(t)=& f_1(t)+f_2(t),
\end{alignat*}
on $(0,(2\pi)^{-n}V(\Omega)]$, where
\begin{alignat*}{1}
f_1(t)&=\xi\omega_n {(nA)^{1+\frac{4}{n}}}t^{-\frac{4}{n}}-\omega_n \frac{4S_{m+2}(nA)^{\frac{n-m+4}{n}}}{(n+4)\rho^{m}}t^{\frac{-2}{n}},\\
\end{alignat*}
and
\begin{alignat*}{1}
f_2(t)&=(1-\xi)\omega_n {(nA)^{1+\frac{4}{n}}}t^{-\frac{4}{n}}+\omega_n \frac{4mS_{m+2}(nA)^{\frac{n-m+3}{n}}}{(n+4)(m+2)\rho^{m+1}}\psi(0)^{\frac{2n-2}{n}}\\
&=(1-\xi)\omega_n {(nA)^{1+\frac{4}{n}}}t^{-\frac{4}{n}}+\omega_n \frac{4mS_{m+2}(nA)^{\frac{n-m+3}{n}}}{(n+4)(m+2)\rho^{m+1}},
\end{alignat*}
for $0<\xi\leq1$.
Then
\begin{alignat*}{1}
\frac{nf_1'(t)}{4\omega_n(nA)^{\frac{4}{n}}}=-\xi(nA)t^{-\frac{n+4}{n}}+\frac{2S_{n+2}}{(n+4)\rho}t^{-\frac{n+2}{n}}.
\end{alignat*}
When
\begin{alignat*}{1}
\frac{2\sqrt{2}S_3}{5k}\leq \xi \leq 1,
\end{alignat*}
we prove that $f(t)$  decreases on $(0,(2\pi)^{-n}V(\Omega)]$ by using
\begin{alignat*}{1}
\frac{\omega_n^{\frac{4}{n}}}{(2\pi)^2} \leq & \frac{1}{2},
\end{alignat*}
and
\begin{alignat*}{1}
\rho\geq (2\pi)^{-n} &\omega_n^{-\frac{1}{n}}V(\Omega)^{\frac{n+1}{n}}.
\end{alignat*}

Therefore,  if
\begin{alignat*}{1}
\frac{2\sqrt{2}S_3}{5} \leq k,
\end{alignat*}
we get
\begin{alignat}{1}
\begin{split}
\sum_{i=1}^k\Gamma_j \geq & \omega_n {(nA)^{1+\frac{4}{n}}}\alpha^{-\frac{4}{n}}-\omega_n \frac{4S_{m+2}(nA)^{\frac{n-m+4}{n}}}{(n+4)\rho^{m}}\alpha^{\frac{mn+m-4}{n}}\\
&+\omega_n \frac{4mS_{m+2}(nA)^{\frac{n-m+3}{n}}}{(n+4)(m+2)\rho^{m+1}}\alpha^{\frac{(m+1)n+m-3}{n}},
\end{split}
\end{alignat}
where
\begin{alignat*}{1}
\alpha=\frac{V(\Omega)}{(2\pi)^n},
\end{alignat*}
and
\begin{alignat*}{1}
\rho&=\frac{V(\Omega)^{\frac{n+1}{n}}}{(2\pi)^n\omega_n^{\frac{1}{n}}}.
\end{alignat*}
Noting that $A=\frac{k}{n\omega_n}$, we obtain the following inequality
\begin{alignat}{1}
\begin{split}
\sum_{i=1}^k\Gamma_j \geq & \omega_n^{-\frac{4}{n}} {\alpha^{-\frac{4}{n}}k^{1+\frac{4}{n}}}
-\omega_n^{\frac{m-4}{n}} \frac{4S_{m+2}}{(n+4)\rho^{m}}\alpha^{\frac{mn+m-4}{n}}k^{\frac{n-m+4}{n}}\\
&+\omega_n^{\frac{m-3}{n}} \frac{4mS_{m+2}}{(n+4)(m+2)\rho^{m+1}}\alpha^{\frac{(m+1)n+m-3}{n}}k^{\frac{n-m+3}{n}}.
\end{split}
\end{alignat}

When $m\geq 2$,  $F(t)$ is defined  by
\begin{alignat*}{1}
F(t)=& F_1(t)+F_2(t)
\end{alignat*}
for $t\in(0,(2\pi)^{-n}V(\Omega)]$, where
\begin{alignat*}{1}
F_1(t)&=\omega_n {(nA)^{1+\frac{4}{n}}}t^{-\frac{4}{n}}+c_3 \omega_n \frac{4mS_{m+2}(nA)^{\frac{n-m+3}{n}}}{(n+4)(m+2)\rho^{m+1}}t^{\frac{(m+1)n+m-3}{n}},\\
\end{alignat*}
for $0<c_3\leq1$ and
\begin{alignat*}{1}
F_2(t)&=-\omega_n \frac{4S_{m+2}(nA)^{\frac{n-m+4}{n}}}{(n+4)\rho^{m}}t^{\frac{mn+m-4}{n}}.
\end{alignat*}
This implies
\begin{alignat*}{1}
\frac{F_1'(t)}{4\omega_n}=&-\frac{1}{n}(nA)^{1+\frac{4}{n}}t^{-\frac{n+4}{n}}\\
&+c_3\frac{(m+1)n+m-3}{n}\frac{mS_{m+2}(nA)^{\frac{n-m+3}{n}}}{(n+2)(m+2)\rho^{m+1}}t^{\frac{mn+m-3}{n}}.
\end{alignat*}
So, if
\begin{alignat*}{1}
c_3\leq \frac{2^{\frac{m+1}{2}}(n+2)(m+2)}{S_{m+2}[(m+1)n+m-3]}k^{\frac{m+1}{n}},
\end{alignat*}
we obtain that $F(t)$  decreases on $(0,(2\pi)^{-n}V(\Omega)]$, which yields that
\begin{alignat}{1}
\begin{split}
\sum_{i=1}^k\Gamma_j \geq & \omega_n^{-\frac{4}{n}} {\alpha^{-\frac{4}{n}}k^{1+\frac{4}{n}}}
-\omega_n^{\frac{m-4}{n}} \frac{4S_{m+2}}{(n+4)\rho^{m}}\alpha^{\frac{mn+m-4}{n}}k^{\frac{n-m+4}{n}}\\
&+c_3\omega_n^{\frac{m-3}{n}} \frac{4mS_{m+2}}{(n+4)(m+2)\rho^{m+1}}\alpha^{\frac{(m+1)n+m-3}{n}}k^{\frac{n-m+3}{n}},
\end{split}
\end{alignat}
where
\begin{alignat*}{1}
c_3\leq \min \left\{ 1, \frac{2^{\frac{m+1}{2}}(n+2)(m+2)}{S_{m+2}[(m+1)n+m-3]}k^{\frac{m+1}{n}} \right\}.
\end{alignat*}
\cvd
\end{Proofp}

For higher eigenvalues, we have the following universal lower bounds
\begin{Corollary}
For any bounded domain $\Omega\subseteq R^n$, $n\geq m\geq 1$ and any $ k \geq 1$ we have\\

(1) When $n=1$ and
\begin{alignat*}{1}
\frac{2\sqrt{2}S_3}{5} \leq k,
\end{alignat*}
we have
\begin{alignat}{1}
\begin{split}
\Gamma_k \geq & \omega_n^{-\frac{4}{n}} {\alpha^{-\frac{4}{n}}k^{\frac{4}{n}}}
-\omega_n^{\frac{m-4}{n}} \frac{4S_{m+2}}{(n+4)\rho^{m}}\alpha^{\frac{mn+m-4}{n}}k^{\frac{-m+4}{n}}\\
&+\omega_n^{\frac{m-3}{n}} \frac{4mS_{m+2}}{(n+4)(m+2)\rho^{m+1}}\alpha^{\frac{(m+1)n+m-3}{n}}k^{\frac{-m+3}{n}},
\end{split}
\end{alignat}\\
where $\alpha$, $\rho$ are defined by (\ref{ar}) and
\begin{alignat*}{1}
S_l=(a+1)^l-a^l.
\end{alignat*}
\\
(2) When $m\geq 2$, we have
\begin{alignat}{1}
\begin{split}
\Gamma_k \geq & \omega_n^{-\frac{4}{n}} {\alpha^{-\frac{4}{n}}k^{\frac{4}{n}}}
-\omega_n^{\frac{m-4}{n}} \frac{4S_{m+2}}{(n+4)\rho^{m}}\alpha^{\frac{mn+m-4}{n}}k^{\frac{-m+4}{n}}\\
&+c_3\omega_n^{\frac{m-3}{n}} \frac{4mS_{m+2}}{(n+4)(m+2)\rho^{m+1}}\alpha^{\frac{(m+1)n+m-3}{n}}k^{\frac{-m+3}{n}},
\end{split}
\end{alignat}
where
\begin{alignat*}{1}
c_3\leq \min \left\{ 1, \frac{2^{\frac{m+1}{2}}(n+2)(m+2)}{S_{m+2}[(m+1)n+m-3]}k^{\frac{m+1}{n}} \right\}.
\end{alignat*}\\
\end{Corollary}

\begin{Proofcp}
Applying Lemma \ref{TLT} and using the similar strategy as in the proof of Theorem \ref{C1}, we prove the desired inequality.

\cvd
\end{Proofcp}

\begin{Remark}
Notice that (\ref{LGIE}) is sharp due to \cite{YY}, our estimate is sharp.\\
\end{Remark}

{\bf Acknowledgments}
 This work was supported supported by the National Natural Science Foundation of China, Grant No. 11531012.  The first author would like to thank Professor Kefeng Liu for his continued support, advice and encouragement.



\bibliography{mybibfile}

\end{document}